\documentclass[12pt,reqno]{amsart}
\usepackage{amsmath, amsthm, amssymb, stmaryrd}
\usepackage{hyperref}
\usepackage{enumerate}
\usepackage{url}
\usepackage{color}
\usepackage{comment}

\let\OLDthebibliography\thebibliography
\renewcommand\thebibliography[1]{
  \OLDthebibliography{#1}
  \setlength{\parskip}{0pt}
  \setlength{\itemsep}{0pt plus 0.3ex}
}

\usepackage{mathtools}

\usepackage{multirow, array}
\usepackage{placeins}

\usepackage{caption} 
\advance \topmargin by -\headheight
\advance \topmargin by -\headsep
     
\evensidemargin \oddsidemargin
\newtheorem{theorem}{\bf Theorem}[section]

\numberwithin{equation}{section}

\newcommand*\wrapletters[1]{\wr@pletters#1\@nil}
\def\wr@pletters#1#2\@nil{#1\allowbreak\if&#2&\else\wr@pletters#2\@nil\fi}

\def \rank {\mathrm{rank}}

\def \dim {\mathrm{dim}}

\def \deg {\mathrm{deg}}

\begin{document}

\title{A remark on sets with few distances in $\mathbb{R}^{d}$}

\author{Fedor Petrov}
\address{St. Petersburg State University, St. Petersburg, Russia}
\email{f.v.petrov@spbu.ru}

\author{Cosmin Pohoata}
\address{California Institute of Technology, Pasadena, CA, USA}
\email{apohoata@caltech.edu}

\thanks{}
\date{}
\begin{abstract} 
A celebrated theorem due to Bannai-Bannai-Stanton says that if $A$ is a set of points in $\mathbb{R}^{d}$, which determines $s$ distinct distances, then $$|A| \leq {d+s \choose s}.$$
In this note, we give a new simple proof of this result by combining Sylvester's Law of Inertia for quadratic forms with the proof of the so-called Croot-Lev-Pach Lemma from additive combinatorics.
\end{abstract}
\maketitle

\section{introduction}

Given a positive integer $s$, a finite subset $A$ in a metric space $M$ is called an $s$-distance set in $M$ if there are $s$ positive real numbers $d_{1},\ldots,d_{s}$ such that all the pairwise distances determined by the points in $M$ are among these numbers, and each $d_{i}$ is realized. Upper bounding the size of such sets is a famous problem in combinatorial geometry, with a lot of activity around the various possible variants. See for instance \cite{GY} and the references therein. When $M$ is $\mathbb{R}^{d}$, with the usual Euclidean distance, the classical result in the area is the following result due to Bannai, Bannai and Stanton \cite{BBS} from 1983.

\bigskip

\begin{theorem} \label{BBS}
If $A$ is an $s$-distance subset in $\mathbb{R}^{d}$, then
$$|A| \leq {d+s \choose s}.$$
\end{theorem}

\bigskip

The proof of Theorem 1 from \cite{BBS} builds upon the linear independence argument introduced for this problem by Larman, Rogers and Seidel in \cite{LRS}. In \cite{LRS}, the authors proved that when $s=2$, the inequality $|A| \leq (d+1)(d+4)/2$ follows from the fact that to each $a$ point in $A$ one can associate a polynomial $f_{a} \in \mathbb{R}[x_{1},\ldots,x_{d}]$ such that $\left\{f_{a}, a \in A\right\}$ is a set of linearly independent polynomials over the reals, which also happens to lie in a subspace of $\mathbb{R}[x_{1},\ldots,x_{d}]$ of dimension $(d+1)(d+4)/2$. This argument was later amplified by Blokhuis \cite{Blok} who showed that one can further add a list of $d+1$ other polynomials to $\left\{f_{a}: a \in A\right\}$ and get an even larger list of linearly independent polynomials that lie in the same vector space of dimension $(d+1)(d+4)/2$. This led to $|A| \leq (d+1)(d+4)/2 - (d+1) = {d+2 \choose 2}$, which established the important first case $s=2$ of Theorem \ref{BBS}. This story was successfully generalized by Bannai-Bannai-Stanton in \cite{BBS}, but for larger $s$ the argument to show that one can add a new list of (higher degree) polynomials to the old list and still get a set of linearly independent elements in the same vector space is significantly more technical.

In this paper, we give a new simple proof of Theorem \ref{BBS} via a slightly improved version of the so-called Croot-Lev-Pach Lemma \cite[Lemma 1]{CLP} over the reals, which may be of independent interest. We state this in a general form, which captures the original version of the Croot-Lev-Pach Lemma as well.

\begin{theorem}\label{CLPPP}
Let $V$ be a finite-dimensional
vector space over a field
$\mathbb{F}$ and $A\subset V$
be a finite set. Let $s$ be a nonnegative integer and let
$p(\overrightarrow{x},\overrightarrow{y})$ be a $2\cdot \dim V$-variate polynomial wih coefficients in $\mathbb{F}$ and of degree at most $2s+1$. Consider the
matrix $M_{p,A}$ with rows and columns indexed by 
$A$ and entries $p(\cdot,\cdot)$. 
It corresponds to a (not necessary
symmetric) bilinear form
on $\mathbb{F}^A$ by a formula
$$
\Phi_p(f,g)= \sum_{a,b\in A} p(\overrightarrow{a},\overrightarrow{b}) f(a)g(b),\, \text{for}\, f,g:A\to \mathbb{F},
$$
which in turn defines a quadratic form 
$\Phi_p(f,f)$.
Denote by $\rank(p,A)$ the rank of matrix
$M_{p,A}$; if $\mathbb{F}=\mathbb{R}$
denote also by $r_+(p),r_-(p)$ the inertia
indices of the quadratic form $\Phi_p(f,f)$. Finally,
denote by $\dim_s(A)$ the dimension
of the space of
polynomials of degree at most
$s$ considered as functions on $A$.
Then:

1) $\rank(p,A)\leqslant 2\dim_s(A)$.

2) if $\mathbb{F}=\mathbb{R}$, then $\max\left\{r_+(p,A),r_-(p,A)\right\}\leqslant \dim_s(A)$.
\end{theorem}

In the next section, we will first prove Theorem \ref{CLPPP}, and then we will use it to deduce Theorem \ref{BBS}. We will need only part 2) of the Lemma above, since part 1) is more or less the original Croot-Lev-Pach lemma in disguise (which doesn't help directly), but we will include nonetheeless a quick new proof of part 1) as well since it motivated part 2).

\section{Proof of Theorem \ref{CLPPP}}

\begin{proof} 
Endow the space $\mathbb{F}^A$ with a natural
inner product $\langle f,g\rangle=\sum_{a\in A} f(a)g(a)$.

Consider the space $\Omega\subset \mathbb{F}^A$ of functions $f$
on $A$ satisfying $\langle f, \phi \rangle=0$
for all polynomials $\phi$ of degree at most
$s$. It is easy to see that the dimension of $\Omega$ as a vector space over $\mathbb{F}$ is at least $|A|-\dim_s(A)$. 

The key observation is that 
$\Phi_p (f,g)=0$ whenever
$f,g\in \Omega$. Indeed, for any
monomial $x^\alpha y^\beta$ in the
polynomial $p(\overrightarrow{x},\overrightarrow{y})$ (here $\alpha,\beta$ are multi-indices
with sum of degrees at most $2s+1$)
we have
$$
\sum_{a,b\in A} a^\alpha b^\beta f(a)g(b)=
\left(\sum_{a\in A} a^\alpha f(a)\right)\cdot 
\left(\sum_{b\in B} b^\beta g(b)\right)=0,
$$
since either $\alpha$ or $\beta$
have degree at most $s$ and $f,g$
are choosing from $\Omega$. 

We will now prove both claims of Theorem \ref{CLPPP} by using dimension arguments. 

Indeed, the bilinear form $\Phi_p[\cdot,\cdot]$ on $\mathbb{F}^A$ takes zero
values on $\Omega \times \Omega$, thus
all non-zero entries of 
its matrix in appropriate basis 
(which includes the basis of $\Omega$
and any other $|A|-\dim\ \Omega$
basis vectors)
may be covered by 
$|A|-\operatorname{dim}\Omega$ rows and $|A|-\operatorname{dim}\Omega$
columns. This implies that every minor of $M_{p,A}$ of dimension at least $2(|A|-\operatorname{dim}\Omega)+1$ must vanish. Therefore,
$$\rank(p,A)\leqslant 2(|A|-\operatorname{dim}\Omega)\leqslant 2\dim_s(A).$$
This proves the first claim of Theorem \ref{CLPPP}.

If $\mathbb{F}=\mathbb{R}$,
by Sylvester's Law of Inertia, we may choose a subspace $Y\subset\mathbb{F}^A$ of dimension
$r_+(p,A)$ such that the
quadratic form $\Phi_p(f,f)$ restricted
to $Y$ is positive definite. 
If $f\in Y\cap \Omega$ and $f\ne 0$, we have
$0=\Phi_p(f,f)>0$, which is impossible.
Therefore, $Y\cap \Omega=\{0\}$ and
$\operatorname{dim} Y+\operatorname{dim}\Omega\leqslant |A|$, which yields that $r_+(p,A)=\dim Y\leqslant |A|-\operatorname{dim}\Omega\leqslant \dim_s(A)$. Analogously, we also have that $r_-(p,A)\leqslant \dim_s(A)$. This completes the proof of Theorem \ref{CLPPP}.
\end{proof}

We now deduce Theorem \ref{BBS} from Theorem \ref{CLPPP}. 

If $A$ is an $s$-distance subset in $\mathbb{R}^{d}$ and $S$ is the set of distinct distances it determines, consider the $2d$-variate polynomial $p$ with real coefficients defined by
$$p(\overrightarrow{x},\overrightarrow{y})=\prod_{d\in S}(d^2-\|x-y\|^2).$$ 
The matrix $M_{p,A}$ from Theorem \ref{CLPPP} is then a positive scalar matrix for this polynomial; therefore, $r_{+}(p,A)=|A|$, and so part 2) of Theorem \ref{CLPPP} implies that 
$$|A|=r_+(p,A)\leqslant \dim_s(A)\leqslant \dim_s(\mathbb{R}^d)={s+d\choose d}.$$
This completes the proof of Theorem \ref{BBS}.

\bigskip

\bigskip

\end{document}